\documentclass[reqno, 11pt]{article}
\usepackage{amssymb,dsfont,graphicx,bbm}
\usepackage{amsthm,amsmath}
\usepackage[usenames,dvipsnames]{color}
\usepackage{fullpage}
\usepackage{ifpdf}
\usepackage[pdftex]{thumbpdf}       

\ifpdf
\usepackage[pdftex]{hyperref}
\hypersetup{%
pdftitle={Large deviations for a random sign Lindley recursion},
pdfauthor={M. Vlasiou, Z. Palmowski},
pdfkeywords={},
bookmarks=true,%
bookmarksnumbered=true,
pdfstartview={FitH},
colorlinks=true,
citecolor=Plum,
linkcolor=RedOrange,
urlcolor=RawSienna,
bookmarksopen=true,%
bookmarksopenlevel=1,     
unicode=true,%
breaklinks=true,%
}%
\else
\newcommand\phantomsection\relax
\newcommand{\url}[1]{#1}
\newcommand{\href}[2]{#2}
\usepackage{nohyperref}
\fi
\theoremstyle{plain}              
\newtheorem{theorem}{Theorem}

\newtheorem{lemma}{Lemma}
\newtheorem{proposition}{Proposition}
\theoremstyle{remark}
\newtheorem{remark}{Remark}

\numberwithin{equation}{section}    

\newcommand{\m}[1]{\mathcal{#1}}
\newcommand{\e}{\mathbb{E}}
\newcommand{\p}{\mathbb{P}}
\renewcommand{\d}{\,\mathrm{d}}

\begin{document}
\title{Tail asymptotics for a random sign Lindley recursion}


\maketitle
\begin{center}
\begin{tabular}{l@{\extracolsep{2cm}}l}
Maria Vlasiou &Zbigniew Palmowski\\
Dept.\ of Mathematics \& Computer Science              &Mathematical Institute\\
Eindhoven University of Technology                     &University of Wroc\l aw,\\
P.O.\ Box 513                                           &pl.\ Grunwaldzki 2/4,\\
5600 MB, Eindhoven                                     &50-384 Wroc\l aw,\\
The Netherlands                                        &Poland\\
\href{mailto:m.vlasiou@tue.nl}{m.vlasiou@tue.nl}       &\href{mailto:zpalma@math.uni.wroc.pl}{zpalma@math.uni.wroc.pl}
\end{tabular}
\end{center}

\begin{abstract}
We investigate the tail behaviour of the steady state distribution of a stochastic recursion that generalises Lindley's recursion. This recursion arises in queuing systems with dependent interarrival and service times, and includes alternating service systems and carousel storage systems as special cases. We obtain precise tail asymptotics in three qualitatively different cases, and compare these with existing results for Lindley's recursion and for alternating service systems.


\end{abstract}

\section{Introduction}\label{s:intro}

This paper focuses on large deviations properties of stochastic recursions of the form
\begin{equation}\label{genlinrec}
W_{n+1}= (Y_nW_n+X_n)^+,
\end{equation}
where $(X_n)$ is an i.i.d.\ sequence of generally distributed random variables, and $(Y_n)$ is an i.i.d.\ sequence independent of $(X_n)$ such that
$\p(Y_1=1)=p=1-\p(Y_1=-1)$, $p\in [0,1]$. Assuming there exists a random variable $W$ such that $W_n \stackrel{\m{D}}{\rightarrow} W$, we are interested in the tail behaviour of $W$, i.e.\ the behaviour of $\p(W>x)$ as $x\rightarrow\infty$. Whitt~\cite{whitt90} has a detailed analysis on the existence of $W$.

The stochastic recursion \eqref{genlinrec} has been proposed as a unification of Lindley's recursion (with $p=1$) and of the recursion
\begin{equation}\label{minusrecursion}
W_{n+1}= (X_n-W_n)^+,
\end{equation}
which is obtained by taking $p=0$. Lindley's recursion~\cite{lindley52} is one of the most studied stochastic recursions in applied probability; Asmussen~\cite{asmussen-APQ} and Cohen~\cite{cohen-SSQ} provide a comprehensive overview of its properties. Recursion \eqref{minusrecursion} is not as well known as Lindley's recursion, but occurs naturally in several applications, such as alternating service models and carousel storage systems. This recursion has been the subject of several studies; see for example \cite{park03,vlasiou07a,vlasiou05,vlasiou05b,vlasiou04}. Most of the effort in these studies has been devoted to the derivation of the distribution of $W$ under various assumptions on the distribution of $X$.

An interesting observation from a methodological point of view is that certain cases that are tractable for Lindley's recursion (leading to explicit expressions for the distribution of $W$) do not seem to be tractable for \eqref{minusrecursion} and vice versa. This was one of the motivations of Boxma and Vlasiou~\cite{boxma07} to investigate the distribution of $W$ defined in \eqref{genlinrec}. For the case $X_n=B_n-A_n$, with $B_n$ and $A_n$ independent non-negative random variables, the work in \cite{boxma07} provides explicit results for the distribution of $W$ assuming that either $B_n$ is phase type and $A_n$ is general, or that $B_n$ is a constant and $A_n$ is exponential.

It appears to be a considerable challenge to obtain the distribution of $W$ under general assumptions on the distribution of $X_n$. To increase the understanding of \eqref{genlinrec} it is therefore natural to focus on the tail behaviour of $W$. Our interest in the tail behaviour of $W$ was raised after realising that the tail behaviour for $W$ can be completely different, depending on whether $p$ is 0 or 1. For example, for $p=0$, it is shown in Vlasiou~\cite{vlasiou07a} that
\begin{equation*}
\p(W>x) \sim \e[e^{-\gamma W}] \p(X>x)
\end{equation*}
if $e^X$ is regularly varying with index $-\gamma$. This includes the case $\gamma=0$ (in which the right tail of $X$ is long-tailed), as well as the case where $X$ has a phase-type distribution (leading to $\gamma>0$). This behaviour is fundamentally different from the case $p=1$, where for example under the Cram\'{e}r condition the tail behaves asymptotically as an exponential; see also Korshunov~\cite{korshunov97} for a concise review of the state of the art. This inspired us to investigate what happens for general $p$.

As is the case for Lindley's recursion, i.e.\ for $p=1$, we find that there are essentially three main cases. For each case, we obtain the asymptotic behaviour of $\p(W>x)$. A brief summary of our results is as follows:
\begin{enumerate}
\item We first consider the case where $X$ has a heavy right tail. In this case, we show that the tail of $W$ is, up to a constant, equivalent to the tail of $X$, under the assumption that $p<1$. Our result shows that there is a qualitative difference with  $p=1$. We derive the tail behaviour of $W$ by developing stochastic lower and upper bounds which asymptotically coincide.

\item The second case we consider is where $X$ satisfies a Cram\'er-type condition, leading to light-tailed behaviour of $W$. By conveniently transforming \eqref{genlinrec}, we are able to apply the framework of Goldie~\cite{goldie91} to get the precise asymptotic behaviour of $\p(W>x)$ as $x\rightarrow \infty$. Our results indicate that for this case there is not a phase transition of the form of the tail asymptotics at $p=1$, but at $p=0$.

\item We finally consider the analogue of the so-called intermediate case, where distributions are light-tailed but the Cram\'er-type condition does not hold. Although the framework of Goldie~\cite{goldie91} does not apply, we can modify some of his ideas to obtain the precise asymptotic behaviour of $\p(W>x)$, assuming that the right tail of $X$ is in the so-called ${\mathcal S}(\gamma)$ class; precise assumptions are stated later in the paper. Interestingly, we find that in this case, there is no phase transition at all; the description of the right tail of $W$ found for $p\in (0,1)$ also holds for the extreme cases $p=0$ and $p=1$.

We believe that the method we develop to deal with the intermediate case is interesting in itself and can also be applied to other stochastic recursions.
\end{enumerate}

This paper is organized as follows. Notation is introduced in Section~\ref{s:notation}. Section~\ref{s:heavy} focuses on the case in which $X$ has a heavy right tail. The Cram\'er case is investigated in Section~\ref{s:light}. The intermediate case is developed in Section~\ref{s:mid}, with which we conclude.

\section{Notation and preliminaries}\label{s:notation}
Throughout this paper, $(X_n)$ and $(Y_n)$ are two mutually independent doubly-infinite i.i.d.\ sequences of random variables as introduced before. We often take generic independent copies $W$ and $X$ of the random variables $W_n$ and $X_n$, as well as of other random variables. Specifically, $W$ is a random variable such that $W \stackrel{\m{D}}{=} (X+YW)^+$, where the random variables $X$, $Y$ and $W$ appearing on the right-hand side are independent. Let $N$ be a random variable such that $\p(N=k)=(1-p)p^{k}$, $k\geqslant 0$ and define $K=N+1$. Loosely speaking, we will use $N$ to count the number of times where $Y=1$ before the event $Y=-1$.

Let $T_n = \sup_{0\leqslant i \leqslant n} S_i$, with $S_0=0$, and $S_i=X_1+\cdots + X_i$ for $i\geqslant 1$. Define the random variable $U$ as $U_i = X_i$ if $Y_i=1$ and $U_i=-\infty$ if $Y_i=-1$. Then
\begin{equation}\label{tnu}
T_N \stackrel{\m{D}}{=} \sup_{n\geqslant 0} [U_0+\cdots +U_n].
\end{equation}
Under the assumption that $\p(X_n<0)>0$ and $\p(Y_1=1)=p=1-\p(Y_1=-1)$, $p\in [0,1)$, which will be made throughout this paper (although we will occasionally compare our results with existing ones for $p=1$), it follows from results in Boxma and Vlasiou~\cite{boxma07} and Whitt~\cite{whitt90} that there exists a stationary sequence $(W_n)$ that is driven by \eqref{genlinrec}; in particular $(W_n)$ is regenerative with finite mean cycle length. As a final point, we use the notational convention $f(x) \sim g(x)$ to denote that $f(x)/g(x)\rightarrow 1$ as $x\rightarrow \infty$.

\section{The heavy-tailed case}\label{s:heavy}
The goal of this section is to obtain the tail behaviour of $W$ assuming that the right tail of $X$ is subexponential. Namely, we assume that the right tail of $X$ is long-tailed. That is, for fixed $y$ it satisfies the following two relations (recall that $X$ may not be positive):
\begin{equation*}
\p(X_1>x) \sim \p(X_1>x+y).
\end{equation*}
and
\begin{equation*}
\p(X_1+X_2>x) \sim 2 \p(X_1>x).
\end{equation*}
We refer to Embrechts {\em et al.}~\cite{embrechts-MEE} for a detailed treatment of subexponential distributions.

The idea of the proof for this case is to first derive stochastic bounds of $W$ in terms of $T_K$ (cf.\ Lemma~\ref{lemht1} below), and then to derive the tail behaviour of $T_K$ (Lemma~ \ref{lemht2}) to obtain the tail behaviour of $W$.

\begin{lemma}\label{lemht1}
It holds that $W\leqslant T_K$ and $W \geqslant T_K - W'$, with both inequalities in distribution, and with $W'$ an independent copy of $W$, independent of $T_K$.
\end{lemma}

\begin{proof}
Consider a stationary version of \eqref{genlinrec}, so that $W \stackrel{\m{D}}{=} W_0$. Note that we can interpret $N$ as
\[
N = \min \{ k \geqslant 0: Y_{-k-1}=-1\},
\]
keeping in mind that we look at events before time zero. Write
\begin{equation*}
\p(W_0>x) = \sum_{n=0}^\infty \p(W_0>x \mid N=n) \p(N=n).
\end{equation*}
The crucial observation is now that the sequence $(W_n)$ behaves like the standard Lindley's recursion between time $-N+1$ and $0$. Since $N$ is a reversed stopping time, it is independent of all events occurring before time $-N-1$. In particular, $N$ is the first time in the past where for $N=k$, $Y_{-k-1}=-1$. That is, $W_{-N} = (X_{-N-1}-W_{-N-1})^+$, and (by stationarity) $W_{-N-1} \stackrel{\m{D}}{=} W$, since $W_{-N-1}$ is determined only by events before time $-N-1$. If we set $W_{k+1}^L = (W_k^L+X_k)^+$, then
\begin{equation*}
\p(W_0>x \mid N=n) = \p(W_0^L> x \mid W_{-n}^L = (X_{-n-1} - W_{-n-1})^+).
\end{equation*}
From this, we see that
\begin{equation*}
\p(W>x) = \sum_{n=0}^\infty \p(W_n^L>x \mid W_0^L = (X-W)^+) \p(N=n).
\end{equation*}
Iterating Lindley's recursion and rearranging indices, we obtain the property
\begin{equation*}
(W_{n}^L \mid W_0^L = (X-W)^+) \stackrel{\m{D}}{=} \max\{0,X_0,\ldots, X_0+\cdots+X_{n-1},X_0+\cdots+X_{n}-W\}.
\end{equation*}
Combining the last two equations, we obtain the bounds
\begin{equation*}
\p(W>x) \leqslant \sum_{n=0}^\infty \p(T_{n+1}>x) \p(N=n),
\end{equation*}
\begin{equation*}
\p(W>x) \geqslant \sum_{n=0}^\infty \p(T_{n+1}-W>x) \p(N=n).
\end{equation*}
The proof follows by noting that $K=N+1$.
\end{proof}

Lemma \ref{lemht1} suggests that the tail behaviour of $W$ is related to the tail behaviour of $T_K$. The tail behaviour of the latter random variable is derived in the next lemma.
\begin{remark}\label{rem:new}
  Note that this result holds without having to make any assumptions on the distribution of $X$. Also, notice that the lemma leads to the alternative lower bound $W\geqslant T_N$.
\end{remark}
Note that
\begin{lemma}\label{lemht2}
If $X$ is subexponential and $p\in [0,1)$, then
\[
\p(T_K>x) \sim \frac 1{1-p} \p(X>x).
\]
\end{lemma}

We omit the proof as the above lemma follows from the fact that the random variable $K$ is independent of the sequence $(X_n)$ and has a light-tailed distribution. In the proof, the sequence of truncated stopping times, and then dominated convergence theorem can be used. For details see Foss and Zachary~\cite{foss03} and Embrechts {\em et al.}~\cite[Lemma 1.3.5]{embrechts-MEE}.

We can now formulate the main result of this section.
\begin{theorem}\label{ht}
If $X$ is subexponential and  $p\in [0,1)$, then
\[
\p(W>x) \sim \frac 1{1-p} \p(X>x).
\]
\end{theorem}

\begin{proof}
Since $T_K$ and $X$ are tail equivalent, and since $X$ is subexponential, $T_K$ is subexponential as well. This implies, in particular, that $T_K$ is long-tailed. This implies in turn that $\p(T_K> x+W') \sim \p(T_K>x)$; see Pitman~\cite{pitman80}. The proof is now completed by invoking Lemma \ref{lemht1}.
\end{proof}

It is interesting to compare Theorem \ref{ht} with existing results for $p=0$ and $p=1$. Theorem~\ref{ht} is consistent with the result $\p(W>x) \sim \p(X>x)$, which holds for $p=0$ and is shown in Vlasiou~\cite{vlasiou07a}, under the assumption that the right tail of $X$ is long-tailed.

As can be expected from the constant $1/(1-p)$, a discontinuity in the asymptotics for $W$ occurs at $p=1$. In this case, it is well known that the asymptotics are of the form $\int_x^\infty \p(X>u) \d u$, which decreases to 0 at a slower rate than $\p(X>x)$; see for example \cite{korshunov97, veraverbeke77, Zachary04} for precise statements.

\section{The Cram\'er case}\label{s:light}

The stochastic bounds of $W$ derived in Lemma~\ref{lemht1} only yield precise asymptotics if $W$ itself is long-tailed. If $X$ has a light right tail, i.e.\ if $\e[e^{\epsilon X}]<\infty$ for some $\epsilon>0$, then $T_K$ satisfies a similar property, implying (by the first part of Lemma~\ref{lemht1}) that $W$ has a light tail as well, which rules out that $W$ is long-tailed. Therefore, we need a different approach to obtain the precise tail asymptotics of $W$. The idea in this section is to relate our recursion to the class of stochastic recursions investigated by Goldie~\cite{goldie91}.

Let $B_n=1$ with probability $p$ and let $B_n=0$ otherwise, where for all $n$, the random variables $B_n$ are independent of each other and of everything else. Define the following three random variables $M_n = B_n e^{X_n}$, $Q_n=e^{X_n}$ and $R_n=e^{W_n}$ and observe that $(M_n,Q_n) \stackrel{\m{D}}{=} (e^{U_n}, e^{X_n})$.

With the obvious notation, we have that
\begin{equation}\label{smarteq}
R \stackrel{\m{D}}{=} \max\{ 1, Q/R, MR\},
\end{equation}
where $Q$, $M$ and $R$ on the right-hand side are independent.
Note that $Q\geqslant M$ a.s. 
We can now obtain the tail behaviour of $R$ by applying Theorem 2.3 of Goldie~\cite{goldie91}. To meet the conditions of Goldie's result, we assume that the distribution of $X$ is non-lattice, and that there exists a solution $\kappa>0$ of $\e[M^\kappa]=1$ satisfying $\e[X e^{\kappa X}]<\infty$, or equivalently
\begin{equation}\label{cramercondition}
\e[e^{\kappa X}]=\frac{1}{p}, \mbox{ such that } m=\e[X e^{\kappa X}]<\infty.
\end{equation}

\begin{theorem}\label{th:cramer}
Under 
condition \eqref{cramercondition} we have that,
\begin{equation*}
\p(R>x) \sim Cx^{-\kappa}, \mbox{ and } \p(W>x) \sim C e^{-\kappa x}
\end{equation*}
with
\begin{equation*}
C= \frac{1}{m} \int_0^\infty [\p(R>t) - \p(MR>t)]t^{\kappa -1} \d t.
\end{equation*}
\end{theorem}

\begin{proof}
The result follows from Theorem 2.3 of Goldie~\cite{goldie91} after we establish that
\begin{equation}\label{conditiongoldie}
\int_0^\infty |\p(R>t) - \p(MR>t)|t^{\kappa -1} \d t<\infty.
\end{equation}
The proof is therefore devoted to verifying \eqref{conditiongoldie}. From \eqref{smarteq} it is clear that $R$ is stochastically larger than $MR$, so we can remove the absolute values in \eqref{conditiongoldie}. Note that
\[
\p(R>t)-\p(MR>t) = \p (\max\{ 1, Q/R, MR\} >t) - \p(MR>t).
\]
Thus, for $t>1$,
\[
\p(R>t)-\p(MR>t) = \p (Q/R >t; MR\leqslant t).
\]
Since $R\geqslant 1$ a.s., this is bounded from above by $\p(Q>t)$.
Thus,
\begin{align}
\nonumber \int_0^\infty |\p(R>t) - \p(MR>t)|t^{\kappa -1} \d t &= \int_0^\infty (\p(R>t) - \p(MR>t))t^{\kappa -1} \d t\\
\nonumber &\leqslant \int_0^1 t^{\kappa-1} \d t + \int_1^\infty t^{\kappa -1} \p(Q>t) \d t \\
\nonumber &\leqslant \frac 1\kappa + \int_0^\infty t^{\kappa -1} \p(Q>t) \d t\\
\label{eq:PRintegr} &= \frac 1\kappa \left(1+  \e[Q^\kappa]\right).
\end{align}
Since $\kappa>0$ and $\e[Q^\kappa]=1/p<\infty$, we conclude that \eqref{conditiongoldie} indeed holds.
\end{proof}

The constant $C$ can be rewritten as follows:
\begin{proposition}\label{prop:constant_cramer}
\begin{equation}\label{Cexpression}
C = \frac {1-p}{m\kappa} + \frac{1-p}m \int_0^\infty \p(X-W>s) e^{\kappa s} \d s + \frac{p}{m} \int_{-\infty}^0 e^{\kappa s} \p(X+W\leqslant s) \d s.
\end{equation}
\end{proposition}

\begin{proof}
Since $R \stackrel{\m{D}}{=} \max \{1,Q/R, MR\}$, we can write
\begin{align*}
\p(R>t)-\p(MR>t) &= \p(\max \{1,Q/R, MR\} > t) -\p(MR>t)\\
&= \p(\max \{1,Q/R\} > t; MR \leqslant t).
\end{align*}
Observe that
\begin{equation}\label{changeofvariableC}
\int_0^\infty t^{\kappa-1} \p(\max \{1,Q/R\} > t; MR \leqslant t) \d t = \int_{-\infty}^\infty e^{\kappa s} \p(\max \{1,Q/R\} > e^s; MR \leqslant e^s) \d s.
\end{equation}

Let $(U, X)$ be a copy of $(U_n,X_n)$; it is useful to recall that $U=X$ with probability $p$ and $U=-\infty$ with probability $1-p$. Since $(Q,M) \stackrel{\m{D}}{=} (e^X, e^U)$, and $W=\log R$, it follows that
\begin{multline*}
\p(\max \{1,Q/R\} > e^s; MR \leqslant e^s) = \p(\max \{0,X-W\} > s; U+W\leqslant s) \\
= (1-p) \p(\max \{0,X-W\} > s) +p \p(\max \{0,X-W\} > s; X+W\leqslant s).
\end{multline*}
Equation \eqref{Cexpression} can now be derived by inserting the above expression in \eqref{changeofvariableC}, distinguishing between positive and negative values of $s$, and some further simplifications.
\end{proof}

Although this provides an expression for the pre-factor $C$, this expression is not explicit as it depends on the entire distribution of $W$. It is therefore interesting to obtain bounds for $C$. From the proof of Theorem~\ref{th:cramer}, it is clear that $C\leqslant \frac{1}{m\kappa} \frac{1+p}{p}$; see also \eqref{eq:PRintegr}. In addition, since $ W\geqslant T_N$ (see Remark~\ref{rem:new}), it is possible to obtain a lower bound for $C$ by deriving the tail behaviour of $T_N$. Specifically, it follows from the representation \eqref{tnu} and a version of the Cram\'er-Lundberg theorem in case the summands of the random walk are equal to $-\infty$ with positive probability that there exists a constant $C_T$ such that
\begin{equation*}
\p(T_N>x) \sim C_T e^{-\kappa x}.
\end{equation*}
This fact actually follows from Theorem 2.3 of Goldie~\cite{goldie91} as well, but can also be proven along the same lines as the standard proof, by mimicking for example the proof of Theorem XIII.5.1 of Asmussen~\cite{asmussen-APQ}. Since $W\geqslant T_N$, we see that $C\geqslant C_T$. Alternative lower and upper bounds may be derived from \eqref{Cexpression}.

Exact computation of $C$ is possible if the exact distribution of $W$ is available. Boxma and Vlasiou~\cite{boxma07} derive expressions for the distribution of $W$ in case $X\stackrel{\m{D}}{=} B-A$, with $B$ a phase-type distribution and $A$ a general distribution. They also obtain the distribution of $W$ in case $B$ is deterministic and $A$ exponential. Computing the exact distribution of $W$ in general seems to be an intractable problem.

As in the previous section, we compare our results with the existing results for $p=0$ and $p=1$. For clarity, write $\kappa=\kappa (p)$ and $C=C(p)$. It is evident that $\kappa(p)$ is continuous at $p=1$ if Equation \eqref{cramercondition} holds for some $p<1$. The constant $C(1)$ can also be shown to be equal to $C_T$, by observing that Theorem XIII.5.1 of Asmussen~\cite{asmussen-APQ} is a special case of Theorem 2.3 of Goldie~\cite{goldie91}. Thus, unlike in the heavy-tailed case, the final asymptotic approximation $C(p) e^{-\kappa (p) x}$ of $\p(W>x)$ is continuous at $p=1$.

Interestingly, it is now the case $p=0$ that is causing some issues. It is shown for the case $p=0$ in Vlasiou~\cite{vlasiou07a} that $$\p(W>x) \sim \e[e^{-\gamma W}] \p(X>x)$$ if $e^X$ is regularly varying of index $-\gamma<0$. If the tail of $X$ is of rapid variation (i.e.\ $\p(X>xy)/\p(X>x)\rightarrow 0$ for fixed $y>1$), then $$\p(W>x) \sim \p(W=0) \p(X>x).$$ Thus, in both cases, the tails of $W$ and $X$ are equivalent up to a constant. From Theorem~\ref{th:cramer} we see that the tail asymptotics for $p>0$ are of a different form.
In particular, if for example $X$ is of rapid variation, then $\e[e^{sX}]<\infty$ for all $s>0$, which implies that Theorem~\ref{th:cramer} holds, i.e.\ the tail of $W$ is exponential for all $p>0$ while it is lighter than any exponential for $p=0$. In this case, as $p\rightarrow 0$,  we have that $\kappa (p)\rightarrow \infty$. In the case where $\e[e^{sX}]=\infty$ for some $s>0$ we distinguish two scenarios.
Let $q=\sup\{s: \e[e^{sX}]<\infty\}$. In the first scenario, the moment generating function is steep, that is $\e[e^{qX}]=\infty$, in which case $\kappa (p)$ converges to $q$. Under the assumptions in  Vlasiou~\cite{vlasiou07a} that $e^X$ is regularly varying of index $-\gamma<0$, we have that $\kappa (p)$ converges to $\gamma$. Note that the asymptotics though might still be of a different form, since $\p(X>x)$ may not have a purely exponential tail.
In the second scenario we have that $\e[e^{qX}]<\infty$. In this case, Theorem~\ref{th:cramer} does not apply if $\e[e^{qX}]$ is less than $1/p$. The study of this case is the subject of the following section.

\section{The intermediate case}\label{s:mid}

In this section we investigate the tail asymptotics when $X$ is light-tailed, but does not satisfy the Cram\'er condition \eqref{cramercondition}. In particular, we assume that $X$ is non-lattice, and a member of the class ${\mathcal S}(\gamma)$ for some $\gamma>0$; that is,
\begin{equation}\label{lgamma}
\p(X>x+y)/\p(X>x) \rightarrow e^{-\gamma y},
\end{equation}
as $x\rightarrow\infty$ and for a fixed $y$, and that
\begin{equation*}
\p(X_1+X_2>x) \sim 2\e[e^{\gamma X}]\p(X>x),
\end{equation*}
where $\e[e^{\gamma X}]<\infty$.
In addition, we assume that
\begin{equation*}
\e[e^{\gamma X}]<1/p,
\end{equation*}
so that the Cram\'er condition \eqref{cramercondition} does not hold (since \eqref{lgamma} implies that $\e[e^{(\gamma+\epsilon) X}]=\infty$ for $\epsilon>0$ and the function $\e[e^{s X}]$ is convex in $s$).

Although the framework of Goldie~\cite{goldie91} does not apply in this case, we are able to modify some of the ideas in that paper to develop an analogue result for the setting of this paper; we believe that our modification is of independent interest.

The main idea is to  derive a useful representation for the distribution of $W$, from which the tail behaviour can be determined.
Define
\[
  g(x)=\p(W>x)-\p(U+W>x),
\]
where $U$ and $W$ are independent, set $V_n = \sum_{i=1}^n U_i$, $V_0=0$ and recall that $S_n=X_1+\cdots + X_n$. The following representation holds.
\begin{lemma}
\begin{equation}\label{sgrep}
\p(W>x)= \p(S_N>x) + \frac{p}{1-p} \p(X+W+S_N\leqslant x; S_N>x) + \p(X-W+S_N>x; S_N\leqslant x).
\end{equation}
\end{lemma}

\begin{proof}
By a telescopic sum argument as in Goldie~\cite[p.\ 144]{goldie91}, we observe that
\begin{align*}
\p(W>x) &=  \sum_{k=1}^n \left( \p(V_{k-1}+ W>x)-\p(V_k+W>x)    \right) + \p(V_n+W >x)\\
&= \sum_{k=1}^n \left(  \p(V_{k-1} + W> x)-\p(V_{k-1} + U+W >x)    \right) + \p(V_n+W >x)\\
&= \sum_{k=0}^{n-1} \int_{-\infty}^\infty g(x-y) \d\p(V_k\leqslant y)+ \p(V_n+W >x).
\end{align*}
Since $V_n\rightarrow -\infty$ a.s.\ as $n\rightarrow\infty$, it follows that
\begin{equation*}
\p(W>x) = \int_{-\infty}^\infty g(x-u) \sum_{n=0}^\infty \d\p(V_n \leqslant u).
\end{equation*}
Note that $g(\infty)=0$, and that the integration range does not include $-\infty$ although $U_i$ does have mass at this point. Moreover, since $\p(V_n \leqslant u) = 1-p^n+p^n\p(S_n\leqslant u)$, we conclude that $\d \p(V_n \leqslant u) = p^n \d\p(S_n\leqslant u)$. Recalling that $\p(N=n)=p^{n}(1-p)$, we obtain
\[
\sum_{n=0}^\infty  p^n \p(S_n \leqslant u) = \frac{1}{1-p} \p(S_N\leqslant u).
\]
Thus,
\begin{equation}\label{intermediateexpression}
\p(W>x)= \frac{1}{1-p} \int_{-\infty}^\infty g(x-u) \d\p(S_N\leqslant u).
\end{equation}
We now simplify the function $g$, using similar arguments as in the proof of Proposition~\ref{prop:constant_cramer} in the previous section.
Note that
\begin{align*}
g(x) &= 1-p+p \p(X+W\leqslant x), \qquad x < 0,\\
g(x) &= (1-p) \p (X-W>x), \qquad x\geqslant 0.
\end{align*}
Inserting this expression for $g$ into \eqref{intermediateexpression}, we obtain
\begin{align*}
\p(W>x) &= \frac 1{1-p} \int_{x}^\infty \left[1-p+p\p(X+W\leqslant x-u)\right]\d\p(S_N\leqslant u)\\
&\quad + \int_{-\infty}^x \p(X-W>x-u) \d\p(S_N\leqslant u)\\
&= \p(S_N>x) + \frac{p}{1-p} \p(X+W+S_N\leqslant x; S_N>x) + \p(X-W+S_N>x; S_N\leqslant x),
\end{align*}
which is identical to \eqref{sgrep}.
\end{proof}

A  crucial second ingredient in obtaining the tail asymptotics is the following useful lemma.

\begin{lemma}\label{sgsn}
If $X \in {\cal S}(\gamma)$ with $\varphi(\gamma)=\e[e^{\gamma X}]<1/p$, then
\begin{align}
\p(S_N>x) &\sim  \frac{(1-p)p}{(1-p\varphi(\gamma))^2} \p(X>x), \\
\p(S_{K}>x) &\sim   \frac{1-p}{(1-p\varphi(\gamma))^2} \p(X>x).
\end{align}
\end{lemma}

\begin{proof}
From, e.g., Cline~\cite[Theorem 1]{cline86}, we have that
\begin{align*}
\p(S_N>x) &\sim  \e[N\varphi(\gamma)^{N-1}] \p(X>x) \\
\p(S_{K}>x) &\sim   \e[K\varphi(\gamma)^{K-1}] \p(X>x).
\end{align*}
Keep in mind that $K=N+1$. The specific constants follow from straightforward computations.
\end{proof}

We are now ready to state and prove the main result of this section.

\begin{theorem}
Let $E_\gamma$ be an exponential random variable of parameter $\gamma $ independent of everything else. Then,
\begin{equation*}
\p(W>x) \sim C_\gamma \p(X>x),
\end{equation*}
with
\begin{equation*}
C_\gamma = \frac{(1-p)p}{(1-p\varphi(\gamma))^2} \left[ \p(X-W+E_\gamma\leqslant 0)+\frac{p}{1-p} \p(X+W+E_\gamma \leqslant 0)\right] +  \frac{1-p}{(1-p\varphi(\gamma))^2} \e[e^{-\gamma W}].
\end{equation*}
\end{theorem}

\begin{proof}
Number the terms in the representation \eqref{sgrep} of $\p(W>x)$ as I, II, III. Lemma \ref{sgsn} yields the tail behaviour of Term I. To obtain the tail behaviour of Terms II and III we make the following useful observation. From Lemma \ref{sgsn} and \eqref{lgamma}, it follows that for fixed $y$,
\begin{equation}\label{expconvergence}
\p(S_N-x>y \mid S_N>x) = \frac{\p(S_N>x+y)}{\p(S_N>x)} \rightarrow e^{-\gamma y} =:\p(E_\gamma >y),
\end{equation}
as $x\rightarrow\infty$. Observe that \eqref{expconvergence} implies
\begin{align*}
\mathrm{II} &= \frac{p}{1-p} \p(X+W+S_N-x\leqslant 0 \mid S_N>x) \p(S_N>x)\\
   &\sim  \frac{p}{1-p} \p(X+W+E_\gamma \leqslant 0) \p(S_N>x).
\end{align*}
For the third term, the argument is similar, but slightly more involved. Write
\begin{align}\label{3ab}
\mathrm{III} &=  \p(X-W+S_N>x; S_N\leqslant x) \nonumber\\
    &= \p(X-W+S_N>x) - \p(X-W+S_N>x; S_N> x).
\end{align}
First, observe that $X+S_N \stackrel{\m{D}}{=} S_{K}$, and observe that the tail of $e^{S_{K}}$ is regularly varying of index $-\gamma$. Since $e^{-W}$ has bounded support, it has finite moments of all orders, so we can apply Breiman's theorem~\cite{denisov07}, as well as the above lemma, to obtain
\begin{equation}\label{3a}
\p(X-W+S_N>x) \sim \e[e^{-\gamma W}] \p(S_K>x).
\end{equation}
To analyse the second term in \eqref{3ab}, observe that
\begin{align*}
\p(X-W+S_N>x; S_N> x) &= \p(X-W+S_N-x>0 \mid S_N> x)\p(S_N> x)\\
&\sim \p(X-W+E_\gamma>0) \p(S_N>x).
\end{align*}
We conclude that
\begin{equation*}
\mathrm{III} \sim \e[e^{-\gamma W}]  \p(S_K>x) -  \p(X-W+E_\gamma>0) \p(S_N>x).
\end{equation*}
Putting everything together, we obtain that
\begin{multline*}
\p(W>x) \sim \left[1+\frac{p}{1-p} \p(X+W+E_\gamma \leqslant 0) -  \p(X-W+E_\gamma>0)\right]\p(S_N>x)+\\ + \e[e^{-\gamma W}]  \p(S_K>x).
\end{multline*}
Simplifying this constant and applying Lemma \ref{sgsn} twice completes the proof.
\end{proof}

Again, we compare our result with the existing results for $p=0$ and $p=1$. For $p=0$, it is shown in Vlasiou~\cite{vlasiou07a} that $\p(W>x) \sim \e[e^{\gamma W}] \p(X>x)$ (note that \eqref{lgamma} guarantees that $e^X$ is regularly varying with index $-\gamma$). This is consistent with the constant $C_\gamma$ defined above, which indeed simplifies to $\e[e^{-\gamma W}]$ when $p=0$.

For $p=1$, it is known (see e.g.\ \cite{bertoin96, korshunov97}) that $\p(W>x) \sim \frac{\e[e^{\gamma W}]}{1-\varphi(\gamma)} \p(X>x)$. This is consistent with our constant $C_\gamma$ specialised to $p=1$, in which case
\[
C_\gamma = \frac {\p(X+W+E_\gamma \leqslant 0) }{(1-\varphi(\gamma))^2}.
\]
In order to show continuity at $p=1$, we need to show that $p=1$ implies  that
$$
\p(X+W+E_\gamma \leqslant 0) = \e[e^{\gamma W}] (1-\varphi(\gamma)).
$$
This can be shown by using the fact that for any non-negative random variable $Y$,
$\e[e^{-\gamma Y}] = \p(Y\leqslant E_\gamma)$, using the identities $e^x+1= e^{x^+}+e^{x^-}$, and $W\stackrel{\m{D}}{=}(W+X)^+$ which holds for $p=1$. To this end, we have that
\begin{align*}
\p(X+W+E_\gamma \leqslant 0) &= \p(E_\gamma \leqslant -(X+W))\\
&= \p(E_\gamma \leqslant -(X+W)^-)\\
&= 1- \p(-(X+W)^- \leqslant E_\gamma)\\
&= 1- \e[e^{\gamma (X+W)^-}]\\
&=  \e[e^{\gamma (X+W)^+}]-  \e[e^{\gamma (X+W)}]\\
&= \e[e^{\gamma W}] -  \e[e^{\gamma X}] \e[e^{\gamma W}]\\
&= \e[e^{\gamma W}] (1-\varphi(\gamma)).
\end{align*}
We conclude that the formula we found for the tail asymptotics in the intermediate case is also valid if $p=0$ or if $p=1$. This contrasts the heavy-tailed case, in which there is a phase transition at $p=1$ (cf.\ Section~\ref{s:heavy}), and the Cram\'er case, where there is a phase transition at $p=0$ (cf.\ Section~\ref{s:light}).

\phantomsection
\addcontentsline{toc}{section}{Acknowledgments}
\section*{Acknowledgments}
The authors are indebted to several colleagues for their help and support during this project. This research originated from a question of Bob Foley. Several helpful comments were provided by Bert Zwart and an anonymous referee. This work is partially supported by the Ministry of Science and Higher Education of Poland under the grant N N2014079 33 (2007-2009).


\phantomsection

\end{document}